\documentstyle[11pt]{article}

\newtheorem{th}{Theorem}[section]

\newtheorem{pro}[th]{Proposition}
\newtheorem{rem}[th]{Remark}

\newcommand{\un}{{\mbox{\tiny $-1$}}}

\newcommand{\uq}{{\sf u}_q\!\left({\cal G}\right)}
\newcommand{\uqp}{{\sf u}_q^+\!\left({\cal G}\right)}
\newcommand{\uqpp}{{\sf u}_q^{++}\!\left({\cal G}\right)}

\newcommand{\vqp}{{\sf v}_q^+\!\left({\cal G}\right)}
\newcommand{\cro}{\underline{\otimes}}

\newcommand{\oti}{{{\mbox{\tiny $\otimes$}}}}

\newcommand{\qu}{Q_{\mbox{\tiny${\cal G}$}}}

\title{Half-quantum groups at roots of unity, path algebras and representation
type}

\author{Claude Cibils}

\begin{document} \date{} \maketitle

\begin{abstract}
Let  ${{\cal G}}$ be a simple Lie algebra of type ${\bf A}$, ${\bf D}$ or ${\bf E}$ and  $q$   a
primitive root of unity of order $n\geq 5$. We show that the finite dimensional half-quantum group
$\uqp$ is of wild representation type, except for ${\cal G} = sl_2$. Moreover, the algebra $\uqp$ is
an admissible quotient of the path algebra of the Cayley graph of the abelian group $\left({\bf
Z}/n{\bf Z}\right)^t$ with respect to the columns of the $t \times t$ Cartan matrix of ${{\cal G}}$.
\end{abstract}

\normalsize

\section{Introduction}

Let ${\cal G}$ be a simple Lie algebra corresponding to a symmetric positive
definite Cartan matrix of type ${\bf A}$, ${\bf D}$ or ${\bf E}$. We consider
in this paper the finite dimensional Hopf algebras  $\uqp$ (see
(\cite{lu2,lu3,lu4,retu1,retu2}) for $q$ a root of unity; they are quotients of
the upper triangular sub-algebra of  the quantum groups  $U_q\left({\cal
G}\right)$ introduced by Drinfel'd and Jimbo \cite{dr1,dr4,ji1,ji2}. Those
structures gives rise to solutions of the Yang-Baxter equation through the
modular category of representations and the universal $R$-matrix acting on it.
They have been  useful in knot theory and 3-manifold invariants, see for
instance  \cite{chpr,ka,tu}.

The representation theory of quantum groups at roots of unity is not fully understood, see
\cite{tu,coka,cokapr,sc}. Actually we will show that  $\uqp$ is of wild representation type except
for ${\cal G}=sl_2$. It is known that finite dimensional algebras are either of finite, tame or wild
representation type, see \cite{dr,cr,ganarosevo}; an algebra is said to be  tame if the isomorphism
classes of indecomposable modules of a fixed dimension are almost all in a finite number of
1-parameter algebraic families. Tame algebras are suitable for reaching a classification of
indecomposable modules, by contrast  the module category of a wild representation type algebra
contains the representation theory of any other finitely generated algebra, see \cite{ga2}. For an
account of the theory, see for instance \cite{dr,pe} or \cite{larasa}. The $sl_2$-case corresponds
to a finite representation type quantum group, see \cite{aqqg}. Suter and Xiao independently
(\cite{su,xi}), have computed the Morita reduction of the entire ${\sf u}_q(sl_2)$, obtaining that
this algebra is of tame representation type, and providing canonical forms for the indecomposable
modules. Notice that the Morita reduction of ${\sf u}_q(sl_2)$ is not anymore a Hopf algebra, hence
the monoidal structure on the canonical forms is unknown.  It is very likely that the entire finite
dimensional quantum groups $\uq$ are of wild representation type if ${\cal G}$ is different from
$sl_2$.

A starting point for our considerations is an easy observation concerning
simple $\uqp$-modules: all of them are one-dimensional. This implies that
$\uqp$ can be presented as a quotient of the path algebra of its canonical
quiver $\qu$ by an ideal of relations $J$  contained in the square of the ideal
generated by the arrows of the quiver and containing some power of it, such a
presentation is called admissible, see \cite{ga1,ga3,gakero}. It insures the
unicity of $\qu$; in turn, simple module of an admissible quotient of a path
algebra are one-dimensional and  in one-to-one correspondence with the vertices
of the quiver. Moreover, the dimension of ${\rm Ext}^1$ between two simple
modules is the number of arrows  between the corresponding vertices.

We describe now the quiver $\qu$, which is a finite oriented graph: let $t$ be
the size of the Cartan matrix $C$ of ${\cal G}$ ($t$ is the number of simple
roots of ${\cal G}$) and let $q$ be a root  of unity of order $n$. The set of
vertices of $\qu$ is the abelian group $G=<K_1,\dots,K_t \mid K_i^n=1,\
K_iK_j=K_jK_i>$. The quiver $\qu$ is the Cayley graph of $G$ with respect to
the set of elements given by the columns of $C$, in other words each vertex
$K^x=K_1^{x_1}K_2^{x_2}\dots K_n^{x_n}$ is the target of $t$ arrows coming from
$\left\{K^{x-a_{\_,j}}\right\}_{j=1,\dots,t}$ where
$\left(a_{\_,j}\right)_{j=1,\dots,t}$ are the columns of $C$ and
$K^{x+a_{\_,j}}$ is the element $K_1^{x_1+a_{1,j}}K_2^{x_2+a_{2,j}}\dots
K_t^{x_t+a_{t,j}}$ of $G$. In case of $sl_2$ the Cartan matrix is $C=(2)$, the
quiver $Q_{_{sl_2}}$ has $n$ vertices $\left\{1,K,\dots,K^{n-1}\right\}$ with
arrows $K^{x-2}\rightarrow K^x$; this quiver is connected for $n$ odd and have
two connected components for $n$ even. Its quantum structure and monoidal
category of representations has been described in \cite{aqqg}, notice that
there are only a finite number of isomorphism classes of indecomposable
 ${\sf u}^+_q(sl_2)$-modules and that their tensor product is given through
Clebsch-Gordan-like formulas.

The Cayley graph of a group is usually  builded using a set of generators of
the group. In our case the columns of $C$ generates $G$ if and only if the
determinant of $C$ is invertible modulo $n$. The quiver $\qu$ is connected
exactly in this case.

The isomorphism $\chi_q : \uqp \longrightarrow k\qu /J$ is given explicitely. The formulas are
provided by structure results about Hopf bimodules  obtained in \cite{acq} and by the quantum
Fourier transform. The latter extends to the level of Hopf bimodules the usual Fourier transform
$\chi : kG \longrightarrow k^G$, i.e. the linear extension of the group isomorphism between a finite
abelian group $G$ and its character group obtained through a choice of generators and  roots of
unity of the field $k$. Notice that $\chi$ is actually a Hopf algebra isomorphism as well as
$\chi_q$ for a natural Hopf structure on $k\qu /J$ that we provide. The way  $\chi_q$ is obtained
guarantees that the formulas gives a well defined isomorphism of Hopf algebras. Alternatively, it is
possible to ignore the origin of the formulas and to verify directly that they provides an
isomorphism of associative algebras. With this point of view, the coalgebra structure on the
admissible path algebra quotient is obtained from $\uqp$ by structure transport through $\chi_q$.

The presentation of $\uqp$ given by the model of paths with comultiplication of $\qu$ is related but
different from the canonical bases of Lusztig (\cite{lu4}). The purpose of this paper is to link
quantum groups at roots of unity with results concerning finite dimensional algebras. C.M. Ringel
has obtained another approach to $\uqp$, by means of indecomposable representations of the quiver
given by the Dynkin diagram associated to $C$ over a finite field, equipped with a Hall algebra like
multiplication, see \cite{ri,gr}.

The representation type of $\uqp$ is obtained in the last section either  through the path
presentation and 2-nilpotent quotients,  or by an induction procedure from a sub-algebra.

We also  use the quiver presentation in order to show that although $\uqp$ is not quasitriangular,
there is an algebra automorphism of a twisted tensor product $\uqp \oti_q \uqp$ transforming  the
comultiplication to its opposite.

{\bf Acknowledgements}: I thank I. Assem and L. Salmer\' on for valuable
exchanges concerning the representation type of some particular algebras.

\section{Quantum paths for half quantum groups}

A half-quantum group is a finite dimensional $k$-Hopf algebra $\uqp$ associated to a symmetric
Cartan matrix $C$ and a root  of unity, the definition is given below. This terminology (see
\cite{crfr}) has been introduced since the ``usual'' quantum groups $\uq$ of Drinfel'd and Jimbo are
easy quotients of the Drinfel'd double of $\uqp$. Half-quantum groups at roots of unity and their
representation theory have direct use in topological quantum field theories, see for instance
\cite{crfr,ku}.

A  symmetric matrix $C=\left(a_{i,j}\right)$ of size $t\times t$ is a Cartan matrix if $a_{i,i}=2$
for $i=1,\dots,t$ and $a_{i,j}\in\{0,-1\}$ for $i\neq j$. Symmetric $t\times t$ Cartan matrices and
graphs with no loops and no multiple edges on $t$ vertices are in natural bijection. Positive
definite Cartan matrices corresponds to graphs of type ${\bf A}$, ${\bf D}$, or ${\bf E}$, see for
instance \cite{bo}. A positive definite symmetric Cartan matrix $C$ provides a semi-simple Lie
algebra ${\cal G}$  with all roots of the same length.

In order to define the corresponding half-quantum group, let $q$ be a root of
unity of order $n\geq 5$ in the ground field. Set $e=n$ if $n$ is odd and
$e=\frac{n}{2}$ if $n$ is even. The associative algebra $\uqp$ has generators
$K_1,\dots,K_t,E_1,\dots,E_t$ subjected to relations that we split in two
sorts:
\begin{itemize}
\item
{\bf I}

\centerline{ $K_i^n=1,\ K_iK_j=K_jK_i,\ K_iE_jK_i^{\un}=q^{a_{i,j}}E_j$}
\item
{\bf II}

\centerline{
 $E_i^e=0,\ E_iE_j=E_jE_i\ \mbox{if}\ a_{i,j}=0,$}

\centerline{$E_i^2E_j-(q+q^\un)E_iE_jE_i+E_jE_i^2=0\ \mbox{if}\ a_{i,j}=-1$.}
\end{itemize}

The comultiplication $\Delta : \uqp \longrightarrow \uqp\otimes\uqp$  is a
coassociative morphism of algebras defined on the generators:

\centerline{$\Delta K_i=K_i\otimes K_i\ \ \ \ \Delta E_i = K_i\otimes E_i +
E_i\otimes 1$.}

The two-sided ideal of relations is preserved by $\Delta$. Lusztig
(\cite{lu2,lu3}) has proved that $\uqp$ is finite dimensional using an action
by algebra automorphisms of the braid group corresponding to $C$, and
constructing root vectors associated to non simple roots.

Our first purpose is to give another presentation of $\uqp$ through the path
algebra of a specific quiver $\qu$. Notice first that  $\uqp$ contains the
group algebra $kG$ of $G=<K_1,\dots,K_t\mid K_iK_j=K_jK_i,\ K_i^n=1>$ as a
sub-Hopf algebra. We assume that the order of the group $\mid\!G\!\mid$ is
invertible in the field $k$. Recall the Fourier transform which is a Hopf
algebra isomorphism
$$\chi : k^G\longrightarrow kG$$
where $k^G=\left\{f:G\rightarrow k\right\}$ is the usual commutative algebra of
functions on a set  provided with the comultiplication given by $\Delta
f(x,y)=f(xy)$ (we identify $k^G\otimes k^G$ and $k^{G\times G}$ since $G$ is
finite), while $kG$ is the vector space of formal linear combinations of
elements of $G$ equipped with the multiplication induced by the group law on
the basis elements and $\Delta s = s\otimes s$ for each $s\in G$. In other
words, $kG$ is the dual Hopf algebra of $k^G$. In this setting we have
$$\chi\left(\delta_{K^c}\right)=\frac{1}{\mid\!G\!\mid}\sum_{K^x\in
G}q^{-cx}K^x$$
$$\mbox{and}\ \ \chi^\un\left(K^a\right)= \sum_{K^x\in G} q^{ax}\delta_{K^x}$$
with the usual convention that $K^a=K_1^{a_1}\dots K_t^{a_t}$ for
$a=(a_1,\dots,a_t)$ and that $\delta_{K^a}$ is the Dirac mass pointing out
$K^a$. Moreover $q^{ax}=q^{a_1x_1+\dots+a_tx_t}$ for $x=(x_1,\dots,x_t)$.

The Fourier transform $\chi$ extends to $\uqp$. In order to make this precise
we consider an intermediate Hopf algebra  $\vqp$ given by the same set of
generators subjected only to relations of type {\bf I}. We assert that $\vqp$
is isomorphic to the path algebra of a quiver $\qu$ (the definition is given
below)  depending only on the Cartan matrix corresponding to ${\cal G}$.

The source for finding $\qu$ and the extension $\chi_q$ of $\chi$ can be
founded in \cite{acq}. Indeed the Hopf algebra $\vqp$ is a tensor algebra over
$kG$ of a $kG$-Hopf bimodule $B$ determined by $E_1,\dots,E_n$. More precisely,
$B$ is a $kG$-bimodule which is also a $kG$-bicomodule with structure maps
compatible with the bimodule structure. The associative tensor algebra $kG\
\oplus \ B\ \oplus \ B\otimes_{kG}B\ \oplus\ \dots$ becomes a Hopf algebra
which is identical to $\vqp$. A theory of $kG$-Hopf bimodules can be developped
(\cite{acq}) and it can be shown that the categories of $kG$ and $k^G$-Hopf
bimodules are equivalent through a Fourier transform functor. Following this
track, the definitions of $\qu$ and $\chi_q$ are dictate. We give them directly
since once the defining formulas are written, they are self-explanatory.

A quiver $Q$ is an oriented graph given by two sets, $Q_0$ the set of vertices,
 $Q_1$ the set of arrows, and two maps $s,t:Q_1\rightarrow Q_0$ providing a
source and a terminus vertex to each arrow. The set $k^{Q_1}$ of functions over the arrows has a
natural bimodule structure over the commutative algebra of functions $k^{Q_0}$: let $\lambda\in
k^{Q_0}$ and $f\in k^{Q_1}$, and the actions be given by $\left(\lambda
f\right)(a)=\lambda\left(t(a)\right)f(a)$ and
$\left(f\lambda\right)(a)=f(a)\lambda\left(s(a)\right)$. In case $Q_0$ and $Q_1$ are finite sets,
the tensor algebra $k^{Q_0}\ \oplus\ k^{Q_1}\ \oplus\ k^{Q_1}\otimes k^{Q_1}\ \oplus\ \dots$ has a
natural basis given by Dirac masses on the paths of the quiver, where a path
$\alpha_n\alpha_{n-1}\dots
\alpha_{1}$ is a sequence of  arrows such that $t(\alpha_i)=s(\alpha_{i+1})$
for each $i$. We denote this associative algebra $k^Q$, it is called the path algebra of $Q$. As a
vector space it coincides with the finite support maps on the set of paths of the quiver. In this
setting commutative algebras corresponds to empty sets of arrows, i.e. by  algebras of maps on sets
(of vertices).

\noindent
{\bf Definition}: Let $C$ be a $t \times t$ symmetric Cartan matrix and $n$ a positive integer
(corresponding to the order of the root  of unity $q$). The quiver $\qu$ has set of vertices the
group $G=<K_1,\dots, K_t \mid K_iK_j=K_jK_i,\ K^n=1>$. There are $t$ arrows ending at each vertex
$t$, which we denote $A\!\left(K^c,1\right),\dots, A\!\left(K^c,t\right)$. They have respectively
directions   $1,\dots ,t$ and sources the vertices $K^{c-a_{\_,1}},
\dots K^{c-a_{\_,t}}$. In other words, $\qu$ is the Cayley graph of $G$ with respect to the sub-set
of $G$ determined by the columns of $C$.

We quote some arithmetic  of the path algebra $k^{\qu}$:

\centerline{$\delta_{K^d}\delta_{A\!\left(K^c,\ i\right)}=
\delta_{A\!\left(K^c,\ i\right)}\ \ \mbox{if $K^d=K^c$ and $0$ otherwise,}$}

\centerline{$\delta_{A\!\left(K^c,\ i\right)}\delta_{K^d}=
\delta_{A\!\left(K^c,\ i\right)}\ \ \mbox{if $K^d=K^{c-a_{i,\_}}$ and $0$
otherwise,}$}

\centerline{$\left\{\delta_{K^c}\right\}$ is a complete set of primitive
orthogonal idempotents.}

\begin{th}
The associative algebra $\vqp$ is isomorphic to the path algebra $k^{\qu}$.
\end{th}
{\bf Proof}: Let $Q$ be the quiver $\qu$ associated to the Cartan matrix $C$.
Since the path algebra is the tensor algebra over $k^{Q_0}$ of the bimodule
$k^{Q_1}$, an algebra map $\phi:k^Q\rightarrow A$ is uniquely and completely
determined by two maps, $\phi_0:k^{Q_0}\rightarrow A$ and
$\phi_1:k^{Q_1}\rightarrow A$ such that the first is an algebra morphism and
the second is compatible with the first, which means that $\phi_1$ is a
$k^{Q_0}$-bimodule map for the bimodule structure of $A$ induced by $\phi_0$.
Recall that the set of vertices of $Q$ is the group $G$. In order to define
$\chi_q:k^Q\longrightarrow \vqp$ consider in  degree  zero  the Fourier
transform
$\left[\chi_q\right]_0:k^{Q_0}=k^G\longrightarrow \vqp$
$$\left[\chi_q\right]_0\left(\delta_{K^c}\right)=
\chi\left(\delta_{K^c}\right)=
\frac{1}{\mid\!G\!\mid}\sum_{K^x\in G}q^{-cx}K^x,$$
and in degree one
$$\left[\chi_q\right]_1\left(\delta_{A\!\left(K^c,\ i\right)}\right)=
\left[\left[\chi_q\right]_0\left(\delta_{K^c}\right)\right]\ E_i.$$
As an an immediate  consequence of the commutation relation between $E_i$ and $K_j$, it is
interesting and easy to verify  that the following equality holds in $\vqp$; it sums up the required
compatibility of $\left[\chi_q\right]_1$ with $\left[\chi_q\right]_0$:
$$\left[\chi_q\right]_0\left(\delta_{K^c}\right)\ \
\left[\chi_q\right]_1\left(\delta_{A\!\left(K^c,\ i\right)}\right)\ \
\left[\chi_q\right]_0\left(\delta_{K^{c-a_{\_,i}}}\right)=
\left[\chi_q\right]_1\left(\delta_{A\!\left(K^c,\ i\right)}\right).$$

Conversely, let $\Psi:\vqp\longrightarrow k^Q$ be defined on the generators by:
$$\Psi \left(K_i\right)=\chi^{\un}\left(K_i\right)=\sum_{K^x\in
G}q^{x_i}\delta_{K^x}$$
$$\Psi\left(E_i\right)=\sum_{K^x\in G}\delta_{A\!\left(K^x,\ i\right)}$$
There is no difficulty for verifying that $\Psi$ factors through the type {\bf I}
relations defining $\vqp$, and that $\chi_q$ and $\Psi$ are inverse one to the
other.

\vskip3mm
\begin{rem}
\rm
In the presentation above of $\vqp$ the root of unity is absent from the
structure constants. Indeed, the Dirac masses on paths of $Q$ make up a
multiplicative basis of $k^Q$: the product of two of them is either the Dirac
mass of the composite path or zero if the paths do not compose.
\end{rem}
\vskip3mm

The structure of Hopf-bimodules over a group algebra obtained in \cite{acq}
provides the path algebra $k^Q$ with a Hopf structure. It  can be recovered by
structure transport through $\chi_q$ and $\Psi$ since we have given explicit
formulas. We summarize this in the following result:

\begin{pro}
Let $q$ be a root  of unity of order $n$, let $C$ be a $t\times t$ symmetric
Cartan matrix and let $Q$ be the correponding quiver, i.e. the Cayley graph of
$G=\left({\bf Z}/n{\bf Z}\right)^t$ with respect to the columns of $C$.

The path algebra $k^Q$ is a Hopf algebra isomorphic to $\vqp$, where the
comultiplication of $k^Q$ is  given by:
$$\Delta\delta_{K^c}=\sum_{K^xK^y=K^c}\delta_{K^x}\oti\delta_{K^y}$$
$$\Delta\delta_{A\!\left(K^c,\ i\right)}= \sum_{K^xK^y=K^c}q^{x_i}\
\delta_{K^x}\oti\delta_{A\!\left(K^y,\ i\right)}\ +\
\sum_{K^xK^y=K^c} \delta_{A\!\left(K^x,\ i\right)}\oti \delta_{K^y}$$
\end{pro}
{\bf Proof}: It is enough to verify that the above formulas correspond through
$\chi_q$ to the given ones for $\vqp$. Consequently they provide a well defined
coassociative algebra morphism $\Delta : k^Q\longrightarrow k^Q\otimes k^Q$.
Observe  that although $k^Q$ is infinite dimensional,  each sum in the formulas
is finite since the group $G$ is finite.

\vskip3mm

A  path of length $m$ in $Q=\qu$ is completely determined by  its target vertex and a sequence of
$m$ directions: $A\!\left(K^c,\ i_mi_{m-1}\dots i_1\right)$ is the path ending  at $K^c$ with
$i_1$-directional first arrow,...,  and $i_m$-directional last arrow. In particular if each arrow of
the path has same direction,  we say that the path {\em has a direction}. We write  $A\!\left(K^c,\
i^m\right)$ such an $i$-directional path. Recall that Dirac masses on paths of $Q$ make up a basis
for $k^Q$.

\begin{th}
\label{paths}
The finite-dimensional half-quantum group $\uqp$ at a root of unity $q$ of
order $n$ is the quotient of the path algebra $k^{\qu}$ by the two sided ideal
$J$ generated by
\begin{itemize}
\item
Dirac masses of each path having a direction and of length $e$ (recall that
$e=n$ if $n$ is odd and $e=\frac{n}{2}$ if $n$ is even).
\item
Commuting square of $\qu$ for each pair of non-related directions, i.e. $\delta_{A\!\left(K^c,\
i,j\right)}=
\delta_{A\!\left(K^c,\ j,i\right)}$ if $a_{i,j}=0$
\item
$\delta_{A\!\left(K^c,\ i,i,j\right)}
-(q+q^\un)\delta_{A\!\left(K^c,\ i,j,i\right)}
+\delta_{A\!\left(K^c,\ j,i,i\right)}$ if $a_{i,j}=-1$.
\end{itemize}
\end{th}
{\bf Proof}: Translate using $\chi_q$, the type {\bf II} relations defining $\uqp$ as a quotient of
$\vqp$ to $k^{\qu}$ .
\vskip3mm

The interest of the above presentation of $\uqp$ is that  $J$ is an {\em
admissible} ideal of the path algebra $k^{\qu}$, which means that $F^m\subset
J\subset F^2$ where $F$ is the two-sided ideal of functions with support
contained in the set  of positive length paths and $m$ is a positive integer.
Indeed, observe that $F^l$ is the two-sided ideal generated by the Dirac masses
of paths of length $l$, hence each generator of $J$ is contained in $F^2$.
Conversely, the image of $F$ in $k^{\qu}/J$ is nilpotent as a consequence of
the Poincar\'e-Birkhof-Witt theorem for $\uqp$, see 5.10 of \cite{lu2}.

An admissible quotient  $k^Q/J$ of a path algebra has a special feature, its Jacobson radical is
$F/J$ as this nilpotent ideal provides a semisimple algebra $\times_{u\in Q_0}k\delta_u$ at the
quotient. Since simple modules for the later are one-dimensional -- the complete list is provided by
$\left\{k\delta_u\right\}_{u\in Q_0}$ -- the same property holds for $k^Q/J$. Of course the radical
$F/J$ has zero action on $k\delta_u$ while $\delta_u$ acts as $1$ and $\delta_v$ acts as $0$ if
$v\neq u$. Moreover it is no difficult to show that the dimension of ${\rm
Ext}^1_{k^Q/J}\left(k\delta_u, k\delta_v\right)$ is the number of arrows from $u$ to $v$. These well
known facts specializes as follows:

\begin{pro}
The irreducible $\uqp$-modules are one-dimensional and in one-to-one
correspondence with $G=<K_1,\dots , K_t \mid K_iK_j=K_jK_i,\   K_i^n=1>$. Let
$S_{K^c}$ and $S_{K^d}$ be  simple modules corresponding to $K^c$ and $K^d$.
Then ${\rm dim}_k{\rm Ext}^1_{\uqp}\left(S_{K^d}, S_{K^c}\right)=1$ if
$K^d=K^{c-a_{\_,i}}$ for some $i$, and $0$ otherwise.
\end{pro}
\vskip3mm
An associative algebra is  {\em connected} if $0$ and $1$ are the only central idempotents.
Otherwise the connected components (or blocks) of the algebra are the sub-algebras in its
decomposition as a product of connected sub-algebras, given by  the complete system of central
primitive idempotents (such a system is unique). In case the algebra is an admissible quotient of a
path algebra of a quiver,  the connected components of the algebra corresponds to the connected
components of the quiver.

\begin{pro}
The algebra $\uqp$ is connected if and only if ${\rm det}\ C$ is invertible in
${\bf Z}/n{\bf Z}$. More precisely, the number of connected components of
$\uqp$ is the cardinality of ${\rm coker}\ C$, considering the Cartan matrix
$C$ as performing an endomorphism of $\left({\bf Z}/n{\bf Z}\right)^t$.
\end{pro}
{\bf Proof}: If $G$ is a group and $S$ is a subset of $G$, the Cayley graph has
set of vertices $G$ and each vertex $g$ is the target of an arrow coming from
$gs^{\un}$ for each $s\in S$. Clearely the number of connected components of
such a graph is the index of the subgroup generated by $S$ in $G$. In our
context $S$ is given by the columns of $C$ as elements of $G=\left({\bf
Z}/n{\bf Z}\right)^t$ and the result follows.
\vskip3mm

We end this section with an observation concerning the comultiplication of $\uqp$, suggested by the
presentation obtained above. In order to simplify the notation, we   denote $\alpha\!\left(K^c,\
i\right)$ the Dirac mass pointing out an arrow $A\!\left(K^c,\ i\right)$. Hence the trivial paths
gives a complete set of primitive orthogonal idempotents, $K^d\alpha\!\left(K^c,\ i\right)=0$ if $K^d\neq
K^c$ and $K^c\alpha\!\left(K^c,\ i\right)=\alpha\!\left(K^c,\ i\right)$. Moreover, $\alpha\!\left(K^c,\
i\right)K^{c-a_{\_,i}}=\alpha\!\left(K^c,\ i\right)$ and $\alpha\!\left(K^c,\ i\right)K^d=0$ if $K^d$ is not
the source vertex of $\alpha\!\left(K^c,\ i\right)$, that is if $K^d\neq K^{c-a_{\_,i}}$.

We will consider a crossed tensor product algebra $\uqp\cro\uqp$. First notice
that the usual tensor product of algebras $\uqp\otimes\uqp$ is also an
admissible quotient of the path algebra of a suitable quiver. Indeed, let
$Q\times Q$ be the quiver with set of vertices $Q_0\times Q_0$ and set of
arrows $Q_1\times Q_0 \coprod Q_0\times Q_1$ together with the obvious source
and terminus maps. The tensor product of algebras $k^Q\otimes k^Q$ can be
identified with the quotient of the path algebra $k^{Q\times Q}$ by the
commuting relations provided by the mixed squares, i.e.
$$(t\alpha, \beta)(\alpha, s\beta)=(\alpha, t\beta)(s\alpha,\beta)$$ for every
couple of arrows $\alpha$ and $\beta$. In order to present $k^Q/J\otimes
k^Q/J$, we just  need  to consider $k^{Q\times Q}$ modulo the preceeding
commuting relations  together with the ``propagation'' of $J$, namely the
two-sided ideal generated by $(J,K^c)$ and $(K^c,J)$ for each vertex $K^c$.

Instead of the commuting relations above, we consider the following relations
provided by each mixed squares:
\small
$$q^{a_{i,j}}
\left(K^c,\  \alpha(K^d,\ j)\right)
\left(( \alpha(K^c,\ i),\ K^{d-a_{\_,j}}\right)
=
\left( \alpha(K^c,\ i),\ K^d\right)
\left(K^{c-a_{\_,i}},\ \alpha(K^d,j)\right)$$
\normalsize
as well as the unchanged quoted propagation of $J$. The admissible quotient of $k^{\qu}\otimes
k^{\qu}$ by those relations presents a crossed tensor product denoted $\uqp\cro\uqp$, which means
that $\uqp\otimes 1$ and $1\otimes \uqp$ are sub-algebras identical to $\uqp$ and $(x\oti 1)(1\oti
y)=(x\oti y)$. However  the reverse tensor product differs, it  is provided by the consequences of
the $q$-commuting mixed square relations above.

\begin{pro}
There exist an algebra automorphism $\phi$ of the crossed product
$\uqp\cro\uqp$ given by $\phi (K^x,\  K^y)=(K^y,\ K^x)$ in degree zero, and by
$$\phi\left(K^x,\ \alpha(K^y,\ i)\right)=
\ q^{-x_i}\ \left(\alpha(K^y,\ i),\ K^x\right)$$
$$\phi \left(\alpha(K^x,\ i),\ K^y\right)=
\ q^{y_i}\ \left(K^y,\ \alpha(K^x,\ i)\right)$$
in degree one, verifying $\phi\Delta=\Delta^{\rm op}$ where $\Delta$ is the
comultiplication of $\uqp$.

\end{pro}

\section{Representation type}

We will study in this Section the representation type of $\uqp$ using two methods which can be
useful in order to produce  families of indecomposable modules.  Actually the first method do not
reach the case $t=3$, i.e. ${\cal G}=sl_3$ or ${\cal G}=sl_2\times sl_2$.

\begin{pro}
\label{nilp}
The algebra $\uqp$ is of wild representation type for $n\geq 5$ and $t\geq 3$.
\end{pro}

{\bf Proof}: We consider the 2-nilpotent quotient of $\uqp$, namely the path
algebra $k^{\qu}$ modulo  $F^2$, which is an algebra having zero square
Jacobson radical. The classification of algebras $k^Q/F^2$ according to their
representation type is known: first consider the separated quiver $Q'$ which
has  double set of vertices $Q_0\times\{0,1\}$ and same set of arrows $Q_1$,
but source and terminus maps given by $s'(a)=(s(a),0)$ and $t'(a)=(t(a),1)$. It
is easy to see that except for simple modules, indecomposable $k^Q/F^2$ and
$k^{Q'}$-modules  corresponds naturally (the second algebra has twice the
number of simple modules of the first one).

Next the representation type of the path algebra $k^{Q'}$ is determined by the
underlying graph of $Q'$ (see \cite{bgp,kac}). If this graph is of Dynkin type
or extended Dynkin (affine) type, the representation type is respectively
finite or tame. Otherwise the representation type is wild. For the Cayley graph
$\qu$ the separated quiver $\qu '$ has at least $3$ arrows leaving each zero
level vertex, which insures that $k^{\qu '}$ is of wild representation type.

\vskip3mm

The second method is based on the usual presentation of $\uqp$, we will show
that the sub-algebra $\uqpp$ generated by $E_1, \dots , E_t$ is already wild,
except for $sl_2$. This is not enough to obtain the same property for the
entire algebra: consider the trivial fact that any finite dimensional algebra
is a sub-algebra of its linear endomorphisms algebra which has only one
indecomposable -- actually simple -- module.  We will see that there exist a
$\uqpp$- sub-bimodule of $\uqp$ complementing $\uqpp$. The next result is due
to  V.M. Bondarenko and Y.A. Drozd, see \cite{bodr} and \cite{larasa}. It
asserts that when a complement as before exists, the sub-algebra provides a
lower bound for the representation type of the entire algebra.

\begin{pro} \label{bound} Let $\Lambda$ be a finite dimensional algebra and let
$A$ be a sub-algebra. Assume there exist an $A$ sub-bimodule $B$ of the
$A$-bimodule $\Lambda$ such that $\Lambda = A \oplus B$. Then if the
representation type of $A$ is wild, the representation type of $\Lambda$ is
also wild. \end{pro}

\vskip3mm
\begin{pro}
The algebra $\uqp$ is of wild representation type for $n\geq 5$ and $t\geq 2$
and of finite representation type for $t=1$, that is for ${\cal G}=sl_2$.
\end{pro}

{\bf Proof}: Consider the quasipolynomial quotient algebra of $\uqpp$  given by
generators $B_1,\dots ,B_t$ subjected to $B_i^e=0$, $B_iB_j=B_jB_i$ if
$a_{i,j}=0$ and $B_iB_j=qB_jB_i$ if $a_{i,j}=-1$ and $i>j$. Since $t\geq 2$ and
$n\geq 5$, we consider further the quotient algebra generated by $C_1$ and
$C_2$, subjected to $C_1^3=C_2^3=0$ and $C_2C_1=C_1C_2$ if $a_{1,2}=0$ or
$C_2C_1=qC_1C_2$ if $a_{1,2}=-1$.
The later has a known wild quotient, namely the algebra (c),  p.238 of
\cite{ri2} on generators $X$, $Y$ and relations $X^2=0$, $XY-\alpha YX=0$,
$Y^2X=Y^3=0$ whith $\alpha\neq 0$. The case of $sl_2$ and $q$ any root of unity
has been studied in \cite{aqqg}, notice however that
${{\sf u}_q^{++}\!\left({sl_2}\right)}
= k[E]/E^n$  is already of finite representation type.

In order to describe a $\uqpp$ sub-bimodule of $\uqp$ complementing $\uqpp$, we
recall the existence of a Poincar\'e-Birkhof-Witt basis. The root vectors
obtained through the action of the braid group corresponding to the Cartan
matrix gives a basis of monomials $K_1^{m_1}\dots K_t^{m_t}E^s$, see for
instance \cite{coka}. A vectorial complement is provided by the linear span of
elements with $m\neq 0$. As a consequence of the commutation relations between
$E_i$ and $K_j$, the later vector space is indeed a $\uqpp$ sub-bimodule.

\tiny

\noindent
D\'epartement  de Math\'ematiques, Universit\'e de Montpellier 2,  F-34095 Montpellier cedex 5

\noindent {\tt cibils@math.univ-montp2.fr}

\vskip1mm
\noindent
Section de math\'ematiques, Universit\'e de Gen\`eve, CP 240, CH-1211 Gen\`eve 24.

\noindent
{\tt Claude.Cibils@math.unige.ch}
\vskip5mm

\noindent January 23, 1996 \end{document}